\documentclass{amsart}

\usepackage{amssymb}

\usepackage[all]{xy}

\newtheorem{thm}{Theorem}

\newtheorem{lem}[thm]{Lemma}
\newtheorem{cor}[thm]{Corollary}

\newtheorem{prop}[thm]{Proposition}

\theoremstyle{definition}
\newtheorem{defn}[thm]{Definition}

\newtheorem{say}[thm]{}
\newtheorem{exmp}[thm]{Example}

\newtheorem{rem}[thm]{Remark}          

\newtheorem*{ack}{Acknowledgments}      
\newtheorem{notation}[thm]{Notation}   
  
\newtheorem{defn-thm}[thm]{Definition--Theorem}  
\newtheorem{defn-lem}[thm]{Definition--Lemma}  

\theoremstyle{remark}

\renewcommand{\o}[0]{{\mathcal O}} 

\renewcommand{\a}[0]{{\mathbb A}}

\newcommand{\p}[0]{{\mathbb P}}

\newcommand{\q}[0]{{\mathbb Q}}
\newcommand{\map}[0]{\dasharrow}
\newcommand{\qtq}[1]{\quad\mbox{#1}\quad}

\newcommand{\gal}[0]{\operatorname{Gal}}

\newcommand{\red}[0]{\operatorname{red}}    
    
\newcommand{\im}[0]{\operatorname{im}}

\newcommand{\aut}[0]{\operatorname{Aut}}

\newcommand{\diff}[0]{\operatorname{Diff}}  
\newcommand{\diffg}[0]{\operatorname{Diff}^*}

\newcommand{\res}[0]{\operatorname{\mathcal R}}

\newcommand{\rdown}[1]{\lfloor{#1}\rfloor}

\newcommand{\onto}[0]{\twoheadrightarrow}

\newcommand{\simq}[0]{\sim_{\q}}

\newcommand{\tsum}[0]{\textstyle{\sum}}

\newcommand{\bir}[0]{\operatorname{Bir}}

\def\into{\DOTSB\lhook\joinrel\to}

\def\loccoh#1.#2.#3.#4.{H^{#1}_{#2}(#3,#4)}

\DeclareMathAlphabet{\mathchanc}{OT1}{pzc}%
                                {m}{it}

\newcommand{\src}[0]{\operatorname{Src}} 
\newcommand{\spr}[0]{\operatorname{Spr}}

\newcommand{\simcb}[0]{\stackrel{cbir}{\sim}}

\begin{document}
\bibliographystyle{amsalpha}


\title{Sources  of log canonical centers}
\author{J\'anos Koll\'ar}

\maketitle

\section{Introduction}

Let $X$ be a smooth variety and $S\subset X$ a smooth hypersurface.
The {\it Poincar\'e residue map} is an isomorphism
$$
\res: \omega_X(S)|_S\cong \omega_S.
$$
In additive form it gives the {\it adjunction formula}
$(K_X+S)|_S\sim K_S$, 
but this variant does not show  that
$\res$ is a {\em canonical} isomorphism.

Its  generalization to log canonical pairs
$(X, S+\Delta)$ has been an important tool in birational geometry;
see, for instance, 
\cite{k-etal, km-book}.
One defines a twisted version of the restriction of $\Delta$
to $S$, called the {\it different} and, for $m>0$ sufficiently divisible,
 one gets a
Poincar\'e residue map
$$
\res^m: \bigl(\omega_X^{[m]}(mS+m\Delta)\bigr)|_S\cong 
\omega_S^{[m]}\bigl(m\diff_S\Delta\bigr),
$$
where the exponent $[m]$ denotes the double dual of the
$m$th tensor power. 
As before, it
 is frequently written as a $\q$-linear equivalence of divisors
$$
\bigl(K_X+S+\Delta)|_S\simq K_S+\diff_S\Delta.
$$
There have been several attempts to extend these formulas
to the case when $S$ is replaced by a higher codimension log canonical center
of a pair $(X,\Delta)$ \cite{kaw1, kaw-adj, k-adj}. 
None of these have been completely successful; the main difficulty
is understanding what kind of object the different should be.

Let $Z\subset X$ be a log canonical center
of a pair $(X,\Delta)$. We can choose a resolution
$f:X'\to X$ such that if we write 
$f^*\bigl(K_{X}+\Delta\bigr)\simq K_{X'}+\Delta'$
then there is a divisor $S\subset X'$ that dominates $Z$
and appears in $\Delta'$ with coefficient 1. The usual adjunction
formula now gives
$$
\bigl(K_{X'}+\Delta'\bigr)|_S\simq K_S+\diff_S(\Delta'-S)=: K_S+\Delta_S.
$$
Note further that $K_{X'}+\Delta'$ is trivial on the fibers of $f$,
hence so is $K_S+\Delta_S$. Thus
$$
f|_S:(S, \Delta_S)\to Z
$$
is a fiber space whose (possibly disconnected)
fibers have (numerically) trivial (log) canonical class.
The aim of previous attempts was to generalize Kodaira's canonical
bundle formula for elliptic surfaces  (cf.\ \cite[Sec.V.12]{bpv}) to this
setting. The difficulty is to make sure that we do not lose
information in the summand that corresponds to the $j$-invariant of
the fibers in the classical case. (For families of elliptic curves
 this could be achieved
by keeping the corresponding variation of Hodge structures as part of our data.)

This suggests that it could be better to view the pair $(S, \Delta_S)$
as the answer to the  problem. However, in general there are
many divisors  $S_j\subset X'$ that satisfy
our requirements and they do not seem to be related to each other
in any nice way.

Our aim is to remedy this problem, essentially by looking at the
smallest possible intersections of the various divisors $S_j$
on a dlt model of $(X,\Delta)$.
There can be many of these models and intersections,
 but they turn out to be birational to
each other and have several unexpectedly nice properties.
These are summarized in the next theorem.
For the rest of this note we work over a field of characteristic 0.

Dlt models,  the different  and
crepant birational equivalence are recalled in
Definitions \ref{doff.defn}--\ref{bir.crep.logstr.say}.

\begin{thm} \label{main.intro.new.thm}
Let  $(X,\Delta)$ be an lc pair,
$Z\subset X$ an lc center and $n:Z^n\to Z$ its normalization.
Let $f:\bigl(X^m, \Delta^m\bigr)\to  (X,\Delta)$ be a dlt model
 and $S\subset X^m$ a minimal (with respect to inclusion)
lc center of $\bigl(X^m, \Delta^m\bigr)$ that dominates $Z$. 
Set $\Delta_S:=\diffg_S\Delta^m$   and $f_S:=f|_S$. 
Let $f^n_S:S\to \tilde Z_S\to Z^n$ denote the Stein factorization.
\begin{enumerate} 
\item (Uniqueness of sources) The  crepant birational  equivalence class of
 $(S,\Delta_S)$ does not depend on the choice of $X^m$ and $S$.
It  is called the {\em source} of $Z$ and
 denoted by 
$\src(Z, X,\Delta)$.
\item (Uniqueness of springs)
 The  isomorphism class of  $\tilde Z_S$ does not depend on the choice
 of $X^m$ and $S$.
It  is called the {\em spring} of $Z$ and denoted by 
$\spr(Z, X,\Delta)$.
\item (Crepant log structure) 
$(S,\Delta_S)$ is dlt,  $K_S+\Delta_S\simq 
f_S^*\bigl(K_{X}+\Delta\bigr)$
and $(S,\Delta_S)$ is klt on the generic fiber of $f_S$. 
\item (Poincar\'e residue map) For $m>0$ sufficiently divisible, there are
 well defined isomorphisms
$$
\begin{array}
{rcl}
 f^*\bigl(\omega_X^{[m]}(m\Delta)\bigr)|_S&\cong & 
\omega_S^{[m]}(m\Delta_S)\qtq{and}\\
 n^*\bigl(\omega_X^{[m]}(m\Delta)|_Z\bigr)&\cong & 
\Bigl(\bigl(f^n_S\bigr)_*\omega_S^{[m]}(m\Delta_S)\Bigr)^{\rm inv}
\end{array}
$$
where  the exponent  {\rm inv} denotes the invariants under the
action of the group of crepant birational self-maps $\bir^c_Z(S,\Delta_S)$.
\item (Galois property) The extension $\tilde Z_S\to Z$ is Galois
and $\bir_Z(S,\Delta_S)\onto \gal\bigl(\tilde Z_S/Z\bigr)$ is surjective.
\item (Adjunction) Assume $\Delta=D+\Delta_1$. Let $n_D:D^n\to D$
be the normalization and $Z_D\subset D^n$ an lc center of
$\bigl(D^n, \diff_{D^n}\Delta_1\bigr)$ such that 
$n_D(Z_D)=Z$. Then 
there is a commutative diagram
$$
\begin{array}{ccc}
\src\bigl(Z_D, D^n, \diff_{D^n}\Delta_1\bigr) &\simcb &
\src\bigl(Z, X, D+\Delta_1\bigr) \\
\downarrow && \downarrow\\
Z_D & \stackrel{n_D}{\to} & Z.
\end{array}
$$
\end{enumerate} 
\end{thm}

Crepant log structures are defined in Section 2.
Theorem \ref{chains.of.lc.centers.thm} shows that
minimal lc centers are birational to each other; this
proves (\ref{main.intro.new.thm}.1) 
and it also establishes (\ref{main.intro.new.thm}.6).
Its consequences for the Poincar\'e residue map
are derived in Section 3.
Sources and springs are formally defined in Section 4
and (\ref{main.intro.new.thm}.5) is proved in Proposition
\ref{min.dlt.source.thm}.

Section 5 contains the  main application, 
Theorems \ref{char.of.slcmod.from.nomr}--\ref{char.of.slcmod.from.nomr.rel}. 
We show that
 normalization gives a one-to-one correspondence:
$$
\left\{
\begin{array}{c}
\mbox{slc  pairs $(X,\Delta)$}\\
\mbox{such that}\\
\mbox{$K_X+\Delta$  is ample}
\end{array}
\right\}
\quad \cong\quad
\left\{
\begin{array}{c}
\mbox{lc  pairs  $\bigl(\bar X, \bar D+\bar\Delta\bigr)$ such that }\\
\mbox{ $K_{\bar X}+ \bar D+\bar\Delta$  is ample plus an }\\
\mbox{involution $\tau$ of
$\bigl(\bar D^n, \diff_{\bar D^n}\bar\Delta\bigr)$}
\end{array}
\right\}.
$$
The papers \cite{odaka, oda-xu}  contain further applications to
$K$-stability and to slc models of deminormal schemes.

Shokurov informed me that his forthcoming paper \cite{shok-source}
 contains another approach to Theorem \ref{main.intro.new.thm}.

\section{Crepant log structures}

\begin{defn}\label{crep.log.str.defn}
 Let $Z$ be a normal variety.
A {\it crepant log structure} on $Z$ is a proper, surjective morphism
$f:(X,\Delta)\to Z$ such that
\begin{enumerate}
\item $f$ has connected fibers,
\item $(X,\Delta)$ is lc and
\item $K_X+\Delta\sim_{f,\q} 0$.
\end{enumerate}

A proper morphism
$f:(X,\Delta)\to Z$ is called a {\it weak crepant log structure} on $Z$
if it satisfies (1) and (3) 
but  $\Delta$
is allowed to be a non-effective sub-boundary.

 Any lc pair $\bigl(Z, \Delta_Z\bigr)$
has a trivial crepant log structure where $(X,\Delta)=\bigl(Z, \Delta_Z\bigr)$.
Conversely, if $f$ is birational then $\bigl(Z, \Delta_Z:=f_*\Delta\bigr)$
is lc.

An irreducible subvariety $W\subset Z$ is a
{\it log canonical center}  or {\it lc center} of
a weak crepant log structure $f:(X,\Delta)\to Z$
 iff it is the image of an
lc center $W_X\subset X$ of $(X,\Delta)$.
A weak crepant log structure has  only finitely many lc centers.

Let $(Z,\Delta_Z)$ be an lc pair and
$f:X\to Z$ a proper, birational morphism.
Write  $K_X+\Delta_X\simq f^*(K_Z+\Delta_Z)$.
Then $f:(X,\Delta_X)\to Z$ is a weak crepant log structure.
The lc centers of $f:(X,\Delta_X)\to Z$ are the same as the
lc centers of  $(Z,\Delta_Z)$.

By Proposition \ref{dlt.mod.exists} we can choose $f$ such that 
 $f:(X,\Delta_X)\to (Z,\Delta_Z)$ is a crepant log structure,
 $X$ is $\q$-factorial and
$(X,\Delta_X)$ is dlt.

Let  $f:(X,\Delta_X)\to Z$ be a dlt crepant log structure
and $Y\subset X$  an lc center.
 Consider the Stein factorization 
$$
f|_Y: Y \stackrel{f_Y}{\longrightarrow} Z_Y\stackrel{\pi}{\longrightarrow} Z
$$ 
and set $\Delta_Y:=\diffg_Y\Delta_X$. 
Then   $\bigl(Y, \Delta_Y\bigr)$ is dlt and
$f_Y:\bigl(Y, \Delta_Y\bigr) \to Z_Y$ is a crepant log structure.
\end{defn}

\begin{defn}[Divisorial log terminal]\label{dlt.defn}
A pair $(X, \sum a_i D_i)$ is called {\it simple normal crossing}
(abbreviated as {\it snc}) 
if $X$ is smooth and for every $p\in X$  one can choose 
an open neighborhood $p\in U$ and local coordinates
$x_i$ such that for every $i$ there is  an index $a(i)$ 
such that $D_i\cap U=(x_{a(i)}=0)$.

As key examples, I  emphasize that the pair
$\bigl(\a^2_k, (x^2=y^2+y^3)\bigr)$ is not snc and 
$\bigl(\a^2_k, (x^2+y^2=0)\bigr)$  is snc 
iff $\sqrt{-1}\in k$. Thus being snc is a
Zariski  local (but not  \'etale local) property.

Given any pair $(X, \Delta)$, there is a largest open subset
$X^{snc}\subset X$ such that 
$\bigl(X^{snc}, \Delta|_{X^{snc}}\bigr)$ is snc. 

A pair $(X, \Delta)$ is called {\it divisorial log terminal}
(abbreviated as {\it dlt}) if the discrepancy
$a(E, X, \Delta)$ is $>-1$ for every divisor whose center is
contained in $X\setminus X^{snc}$. 
\end{defn}

\begin{defn}[Different]\label{doff.defn}
Let $(X,\Delta)$ be a dlt pair and $Y\subset X$  an lc center.
Generalizing the usual notion of the different \cite[Sec.16]{k-etal},
there is a naturally defined $\q$-divisor $\diffg_Y\Delta$,
called the {\it different} of $\Delta$ on $Y$ such that
$$
\bigl(K_X+\Delta\bigr)|_Y\simq K_Y+\diffg_Y\Delta.
$$
The traditional different  \cite[Sec.16]{k-etal} is defined such that
if $Y=D$ is a divisor then
$$
\bigl(K_X+D+\Delta\bigr)|_D\simq K_D+\diff_D\Delta.
$$
Thus, in this case, $\diffg_D(D+\Delta)=\diff_D\Delta$.
This  inductively defines $\diffg_Y\Delta$ whenever
$Y$ is an irreducible component of 
a complete intersection of divisors in $\rdown{\Delta}$.
In the dlt case, this takes care of  every lc center by
\cite[Sec.3.9]{fuj-in-corti};
see also \cite[Sec.4.1]{kk-singbook} for details.
\end{defn}

The following result was  proved by Hacon (and published in
\cite{k-db}).  A simplified proof is in \cite{fujino-ssmmp}. 

\begin{prop}\label{dlt.mod.exists}
 Let $(Z,\Delta_Z)$ be an lc pair. 
Then it has a $\q$-factorial, crepant, dlt model
$p:(X,\Delta_X)\to (Z,\Delta_Z)$.
That is,  $X$ is $\q$-factorial, $(X,\Delta_X)$ is dlt,
$K_{X}+\Delta_X$ is $p$-nef and
$\Delta_X=E+p^{-1}_*\Delta_Z$ where $E$ contains all $p$-exceptional divisors
with multiplicity 1.
\qed
\end{prop}

\begin{say}[Birational weak crepant log structures]
\label{bir.crep.logstr.say}{\ }

Let $f:(X,\Delta)\to Z$ be a weak crepant log structure.
If $f$ factors as $X\stackrel{g}{\to} X'\stackrel{f'}{\to} Z$
where $g$ is birational, 
then $f':(X',\Delta':=g_*\Delta)\to Z$ also 
a weak crepant log structure.
 We say that
$f:(X,\Delta)\to Z$ {\it birationally dominates} $f':(X',\Delta')\to Z$.

Conversely, assume that 
 $f':(X',\Delta')\to Z$  is a weak crepant log structure
and $g:X\to X'$ is a proper   birational morphism.
Write $K_X+\Delta\simq g^*\bigl(K_{X'}+\Delta'\bigr)$.
Then $f:=f'\circ g: (X,\Delta)\to Z$ is also a weak crepant log structure.

By Proposition \ref{dlt.mod.exists} every (weak) 
crepant log structure $f:(X,\Delta)\to Z$
is dominated by another (weak) crepant log structure $f^*:(X^*,\Delta^*)\to Z$
such that $(X^*,\Delta^*)$ is dlt and $\q$-factorial.
If $\Delta$ is effective then we can choose $\Delta^*$ to be effective.

Two weak crepant log structures $f_i:(X_i,\Delta_i)\to Z$ are called
{\it crepant birational} if 
there is a third weak crepant log structure $h:(Y,\Delta_Y)\to Z$
which birationally dominates both of them. Crepant birational 
equivalence is denoted by $\simcb$.


The group of crepant birational self-maps of a weak crepant log structure
$f:(X,\Delta)\to Z$ is denoted by $\bir^c_Z(X,\Delta)$. 
By also allowing $k$-automorphisms, we get the larger group 
 $\bir^c_k(X,\Delta)$.

Let  $f:(X,\Delta)\to Z$ be a weak crepant log structure and
$f':X'\to Z$ a proper morphism. Assume that there is a birational map
$\phi:X\map X'$ such that $f'\circ \phi=f$. By the above, there is a
unique  $\q$-divisor $\Delta'$ such that
 $f':(X',\Delta')\to Z$ is a weak crepant log structure that is birational to
 $f:(X,\Delta)\to Z$. 
If $\phi^{-1}$ has no exceptional divisors, then
$\Delta'=\phi_*\Delta$ and hence
 $\Delta'$ is effective if $\Delta$ is.
 
 Let  $f_i:(X_i,\Delta_i)\to S$ be   weak crepant log structures
and $\phi:X_1\map X_2$ a birational map. Let
$Z_1\subset X_1$ an lc center such that,   at the generic point of $Z_1$,
the pair  $(X_1,\Delta_1)$
 is dlt and $\phi$ is a local
isomorphism. Then
$Z_2:=\phi_*Z_1$ is also an lc center  and
$$
\phi|_{Z_1}:\bigl(Z_1, \diffg_{Z_1}\Delta_1\bigr)\map
\bigl(Z_2, \diffg_{Z_2}\Delta_1\bigr) \qtq{is crepant birational.}
$$
\end{say}


 \begin{thm}\cite{nak-ue, ueno, gon-ab, fuj-gon-cbf} \label{dlt.source.defn.3}  
Let $f:(X,\Delta_X)\to Z$ be a crepant log structure.
 Then: 
\begin{enumerate}
\item The $\bir^c_Z (X,\Delta_X)$ action   on 
$\omega_X^{[m]}(m\Delta_X)$ is finite for every $m\geq 0$.
\item 
If $Z$ is projective
and $K_X+\Delta_X\simq f^*(\mbox{ample $\q$-divisor})$ then the
 $\bir^c_k (X,\Delta_X)$ action  on $Z$ is finite. \qed
\end{enumerate}
\end{thm}

\begin{say}[Minimal dominating lc centers]\label{mincent.bir.inv}
Let  $f:(X,\Delta)\to S$ be a  dlt,  weak crepant log structure.
Let $W\subset S$ be an lc center and
$\{W_i:i\in I(W)\}$  the minimal (with respect to inclusion) 
 lc centers of $(X,\Delta)$ that dominate $W$.
We claim that the set of their crepant birational isomorphism classes
$$
\bigl\{\bigl(W_i, \diffg_{W_i}\Delta\bigr): i \in I(W)\bigr\}
\eqno{(\ref{mincent.bir.inv}.1)}
$$
is a birational invariant of  $f:(X,\Delta)\to S$.

To see this note that by \cite{szabo}
we can assume that $(X,\Delta)$ is   snc. 
Then it is enough to check birational invariance for one
smooth blow up. If we blow up  $V\subset X$ that is not an lc center,
then the set of lc centers is unchanged.

If $V$ is an lc center that is the
complete intersection of say $D_1,\dots, D_r\subset \rdown{\Delta}$,
then we get an exceptional divisor $E_V$ that is a $\p^{r-1}$-bundle
over $V$. Locally on $V$, we get a direct product
$$
\bigl(E_V, \diffg_{E_V}\Delta_{B_VX}\bigr)\cong 
\bigl(V, \diffg_{V}\Delta\bigr)\times
\bigl(\p^{r-1}, (x_1\cdots x_r=0)\bigr),
$$
thus  every minimal lc center of $\bigl(V, \diffg_{V}\Delta\bigr)$
corresponds to  $r$ isomorphic copies of itself 
among the  minimal lc centers of
$\bigl(E_V, \diffg_{E_V}\Delta_{B_VX}\bigr)$, hence
among the  minimal lc centers of $\bigl(B_VX, \Delta_{B_VX}\bigr)$.\qed
\end{say}

Our next aim is to prove that
for crepant log structures,
the  invariant defined in 
(\ref{mincent.bir.inv}.1) consist of a single
birational equivalence class.

\subsection*{$\p^1$-linking of minimal lc centers}{\ }

\begin{defn}[$\p^1$-linking]\label{p1-link.defn} 

A {\it standard $\p^1$-link} is 
 a  dlt, $\q$-factorial, pair $\bigl(X, D_1+D_2+\Delta)$
whose sole lc centers are $D_1, D_2$ (hence $D_1$ and 
$ D_2$ are disjoint)
plus a proper morphism $\pi:X\to S$
such that  $K_X+D_1+D_2+\Delta\sim_{\q,\pi}0$, 
$\pi:D_i\to S$ are both isomorphisms
and every reduced fiber  $\red X_s$ is isomorphic to $\p^1$.

Let $F$ denote a general smooth fiber. Then
$\bigl((K_X+D_1+D_2)\cdot F\bigr)=0$, hence $(\Delta\cdot F)=0$.
That is, $\Delta$ is a vertical divisor, 
the projection gives an isomorphism
$\bigl(D_1, \diff_{D_1}\Delta\bigr)\cong \bigl(D_2, \diff_{D_2}\Delta\bigr)$
and these pairs are klt.

The simplest example of a standard $\p^1$-link is a product
$$
\bigl(S\times \p^1, S\times \{0\}+S\times \{\infty\}+
\Delta_S\times \p^1\bigr)
$$
for some $\q$-divisor $\Delta_S$.

It turns out that
 every standard $\p^1$-link is locally the quotient of a
product. 
To see this note  that $\bigl((D_1-D_2)\cdot F\bigr)=0$, thus every
point $s\in S$ has an open neighborhood $U$ such that
$D_1-D_2\simq 0$ on $\pi^{-1}(U)$. Taking the corresponding cyclic cover
we get another standard $\p^1$-link
$$
\tilde\pi:\bigl(\tilde X_U, \tilde D_1+\tilde D_2+\tilde \Delta)\to \tilde U
$$
where the $\tilde D_i$ are now Cartier divisors and
$\tilde \Delta=\tilde\pi^*\tilde \Delta_U$ for some $\q$-divisor
$\tilde \Delta_U$.  Here 
$\tilde D_1\sim\tilde D_2$, hence 
the linear system $|\tilde D_1,\tilde D_2|$ maps $\tilde X_U$ to $\p^1$.
Together with $\tilde \pi$ this gives
 an isomorphism
$$
\bigl(\tilde U\times \p^1, \tilde U\times \{0\}+\tilde U\times \{\infty\}+
\tilde \Delta_U\times \p^1\bigr)\cong
\bigl(\tilde X_U, \tilde D_1+\tilde D_2+\tilde \Delta).
$$

Let $g:(X,\Delta)\to S$ be  a crepant, dlt  log structure 
and $Z_1, Z_2\subset X$   two lc centers. We say that  $Z_1, Z_2$ 
are {\it directly $\p^1$-linked} if there is an lc center
$W\subset X$ containing the $Z_i$ such that
 $g(W)=g(Z_1)=g(Z_2)$ and
$\bigl(W, \diffg_W\Delta\bigr)$ is 
crepant birational to a
standard $\p^1$-link with $Z_i$ mapping to $D_i$.

We say that  $Z_1, Z_2\subset X$
are {\it  $\p^1$-linked } if there is a sequence of
lc centers $Z'_1,\dots, Z'_m$ such that $Z'_1=Z_1, Z'_m=Z_2$
and $Z'_i$ is  directly $\p^1$-linked to $Z'_{i+1}$
for $i=1,\dots, m-1$ (or $Z_1=Z_2$).
\end{defn}

The following strengthening  of \cite[1.7]{k-db}
 was the reason to introduce the notion of
  $\p^1$-linking.

\begin{thm}\label{chains.of.lc.centers.thm}
Let $k$ be a field  and $S$ essentially of finite 
type over $k$.
 Let $f:(X,\Delta)\to S$ be a proper morphism 
such that $K_X+\Delta\sim_{\q,f}0$ and
 $(X,\Delta)$ is dlt.
Let $s\in S$ be a point 
such that $f^{-1}(s)$ is connected (as a $k(s)$-scheme).
Let
 $Z\subset X$ be minimal (with respect to inclusion)
 among the lc centers  of  $(X,\Delta)$ 
such that  $s\in f(Z)$. Let  $W\subset X$ be an lc center  of  $(X,\Delta)$ 
such that  $s\in f(W)$.

Then there is an  lc center $Z_W\subset W$ such that 
 $Z$ and $Z_W$ are $\p^1$-linked.

In particular, all the minimal (with respect to inclusion)
 lc centers  $Z_i\subset X$ 
such that  $s\in f(Z_i)$ are  $\p^1$-linked to each other.
\end{thm}

Remarks. For the applications it is crucial to understand the case
when $k(s)$ is not algebraically closed. 

Each $\p^1$-linking defines a birational
map $Z\map Z_W$, but different $\p^1$-linkings can give
different birational maps.

\medskip

Proof.    We use induction on $\dim X$ and on $\dim Z$.

Write $\rdown{\Delta}=\sum D_i$. 
By passing to a suitable \'etale neighborhood of $s\in S$ we may assume that
each $D_i\to Y$ has connected fiber over $s$
and every lc center of $(X,\Delta)$ intersects $f^{-1}(s)$.
(We need to do this without changing the residue field
 so that $f^{-1}(s)$ stays connected,
cf.\ \cite[I.4.2]{milne}.)

Assume first that 
 $f^{-1}(s)\cap \sum D_i$ is connected. 
By suitable indexing, we may assume that
$Z\subset D_1$,
$W\subset D_r$ and  $f^{-1}(s)\cap D_i\cap D_{i+1}\neq\emptyset$
for $i=1,\dots, r-1$.

By  induction, we can apply Theorem \ref{chains.of.lc.centers.thm} to
$D_1\to S$ with $Z$ as $Z$ and $D_1\cap D_2$ as $W$.
We get that there is an  lc center $Z_2\subset W$ such that 
 $Z$ and $Z_2$ are $\p^1$-linked. 
As we noted in Definition \ref{p1-link.defn},  
$Z_2$ is also minimal (with respect to inclusion)
 among the lc centers  of  $(X,\Delta)$ 
such that  $s\in f(Z_2)$. 
Note that $Z_2$ is an lc center of
$\bigl(D_1, \diffg_{D_1}(\Delta)\bigr)$.
By adjunction, it is an lc center of $(X,\Delta)$
and  also an lc center of
$\bigl(D_2, \diffg_{D_2}(\Delta)\bigr)$.

Next we apply Theorem \ref{chains.of.lc.centers.thm} to
 $D_2\to S$ with $Z_2$ as $Z$ and 
$D_2\cap D_3$ as $W$, and so on.
At the end we work on $D_r\to S$ with $Z_r$ as $Z$ and 
$W$ as $W$ to get an lc center $Z_W\subset W$ such that 
 $Z$ and $Z_W$ are $\p^1$-linked.
This proves the first claim if  $f^{-1}(s)\cap \sum D_i$ is connected.

If $f^{-1}(s)\cap \sum D_i$ is disconnected, then  
write $\Delta=\sum_{i=1}^m D_i+ \Delta_1$.
We claim that in this case $m=2$ and 
 $D_1, D_2\subset X$ are directly  $\p^1$-linked
(by $W=X$). We may assume that $X$ is $\q$-factorial.

First we show that $\sum D_i$ dominates $S$.
Indeed, consider the exact sequence
$$
0\to \o_X(-\tsum D_i)\to \o_X\to \o_{\sum D_i}\to 0
$$
and its push-forward
$$
\o_S\cong f_*\o_X\to f_*\o_{\sum D_i}\to R^1f_*\o_X(-\tsum D_i).
$$
Since $-\sum D_i\sim_{\q,f} K_X+\Delta_1$, the sheaf
$R^1f_*\o_X(-\sum D_i)$ is torsion free by 
\cite{k-hdi1} (see \cite{k-db} for the extension to the
klt case that we use).
Thus $\o_S\onto f_*\o_{\sum D_i}$ is surjective hence
$ \sum D_i\to S$ has connected fibers, a contradiction.

This implies that $K_X+\Delta_1$ is not $f$-pseudo-effective
and so by \cite[1.3.2]{bchm}
one can run the $(X,\Delta_1)$-MMP over $S$.

  Every step is numerically $K_X+\sum D_i+\Delta_1$-trivial, hence
 $\sum D_i$ is ample on every extremal ray. Therefore a
connected component of $\sum D_i$ can never be contracted by a birational
contraction. By   the
Connectedness Theorem \cite[17.4]{k-etal}, the connected components
of  $\sum D_i$ are unchanged
for birational contractions and flips.
 Thus, at some point, we
  must encounter a Fano contraction 
$p:(X^*,\Delta_1^*)\to V$ 
where $\sum D_i^*$ is $p$-ample.  
So there is an irreducible component, say $D_1^*$ that 
has positive intersection with the contracted ray. Therefore  $D_1^*$ is
  $p$-ample.
By assumption, there is another irreducible component, say $D_2^*$ that is
disjoint from  $D_1^*$.
Let $F_v\subset X^*$ be any fiber that intersects $D_2^*$.
Since  $D_2^*$ is disjoint from  $D_1^*$, we see that $D_2^*$
does not contain $F_v$. Thus $D_2^*$ also 
has positive intersection with the contracted ray, hence  $D_2^*$ is
also  $p$-ample. 

Thus  $D_1^*$ and $D_2^*$ are both relatively ample 
(possibly multi-) sections of $p$ and they are disjoint. 
This is only possible if $p$ has fiber dimension 1,
 the generic fiber is a smooth rational curve and
$D_1^*$ and $D_2^*$ are sections of $p$. 

Since $p$ is an extremal contraction, $R^1p_*\o_{X^*}=0$, which implies that
every fiber of $p$ is a tree of smooth rational curves.
Both $D_1^*$ and $D_2^*$ intersects every fiber in a single point
and they both intersect every contracted curve.
Thus every fiber is irreducible and so
$p:(X^*,\Delta^*)\to V$ is a standard $\p^1$-link
with $D_1^*,D_2^*$ as sections.
As we noted in Definition \ref{p1-link.defn}, the rest of $ \Delta^*$
consists of vertical divisors. Thus any other $D_i^*$  would make
$f^{-1}(s)\cap \sum D_i$ connected. Therefore
 $D_1^*,D_2^*$ are the only lc centers of 
$(X^*,D_1^*+D_2^*+\Delta_1^*)$
and so  $D_1,D_2$ are the only lc centers of 
$(X,\Delta)$. 
As noted at the end of Definition \ref{bir.crep.logstr.say},
 $D_1, D_2\subset X$ are directly  $\p^1$-linked
(by $W=X$).
\qed
\medskip

\begin{cor}\label{lc.cnet.crep.adj.cor} Let 
 $f:(X,\Delta_X)\to S$ be a dlt,   crepant log structure.
Let $Y\subset X$ be an lc center.
Consider the Stein factorization 
$f|_Y: Y \stackrel{f_Y}{\longrightarrow} S_Y\stackrel{\pi}{\longrightarrow} S$ 
and set $\Delta_Y:=\diffg_Y\Delta_X$.  Then
\begin{enumerate}
\item $f_Y:\bigl(Y, \Delta_Y\bigr) \to S_Y$ is a dlt, crepant log structure.
\item Let $W_Y\subset S_Y$ be an lc center of 
$f_Y:\bigl(Y, \Delta_Y\bigr) \to S_Y$. Then
$\pi(W_Y)\subset S$ is an lc center of $f:(X,\Delta_X)\to S$ and
every  minimal lc center of $\bigl(Y, \Delta_Y\bigr)$
dominating $W_Y$ is also a  minimal lc center of $(X,\Delta_X)$
dominating $\pi(W_Y)$.
\item  Let $W\subset S$ be an lc center of $f:(X,\Delta_X)\to S$.
Then every irreducible component of $\pi^{-1}(W)$ is  an lc center of 
$f_Y:\bigl(Y, \Delta_Y\bigr) \to S_Y$.
\end{enumerate}
\end{cor}

Proof. (1) is clear.
To see (2), note that $ W_Y$ is dominated by an lc center $V_Y$ of
$\bigl(Y, \diffg_Y\Delta)$. Thus, by adjunction,
$V_Y$ is also an lc center  of
$(X, \Delta)$, hence $\pi(W_Y)=f(V_Y)$ is an  lc center of $S$.
By Theorem 
\ref{chains.of.lc.centers.thm}, a minimal lc center of $Y$ that dominates
$W_Y$ is also a minimal lc center of $X$ that dominates
$\pi(W_Y)$. Thus $\src\bigl(W_Y, Y, \Delta_Y\bigr)\sim 
\src\bigl(\pi(W_Y),X,\Delta_X\bigr)$.

Finally let $W\subset S$ be an lc center of $f:(X,\Delta_X)\to S$
and $w\in W$ the generic point. 
Let $V_X\subset X$ be a minimal lc center that dominates $W$.
By Theorem \ref{chains.of.lc.centers.thm}, there is an
lc center $V_Y\subset Y$ that is $\p^1$-linked to $V_X$. 
By adjunction, $V_Y$ is also an lc center
of $\bigl(Y, \diffg_Y\Delta)$. Thus
$f_Y(V_Y)\subset S_Y$ is an  lc center of
$f_Y:\bigl(Y, \Delta_Y\bigr) \to S_Y$ and it is also
one of the irreducible components of  $\pi^{-1}(W)$.

In order to get (3), after replacing $S$ by an \'etale neighborhood
of $w$, we may assume that  $Y=\cup Y_j$  such that each 
$f^{-1}(w)\cap Y_j$ is connected.
By the previous argument, each $Y_j$  yields an
lc center $f_{Y_j}(V_{Y_j})\subset S_{Y_j}$ and together
these show that 
every irreducible component of $\pi^{-1}(W)$ is  an lc center of 
$f_Y:\bigl(Y, \Delta_Y\bigr) \to S_Y$.\qed
\medskip

\begin{exmp}\label{difefrnt.p1.links.exmp}
 Fix  $m\ge 3$  and $\epsilon$ a primitive $m$th root of unity.
On $\p^{m-1}$ consider the $\mu_m$-action generated by
$$
\tau_1:(x_0:x_1:\cdots : x_{m-1})\mapsto 
(x_0:\epsilon x_1:\cdots :\epsilon^{m-1} x_{m-1}).
$$
The action moves the divisor
$D_0:=(x_0+x_1+\cdots+x_{m-1}=0)$ into $m$ different divisors
$D_0,\dots, D_{m-1}$. One easily  checks that
$\bigl(\p^{m-1}, D_0+\dots + D_{m-1}\bigr)$  is snc 
(if $\epsilon$ is in our base field) and
has trivial log canonical class.

Let $A$ be an abelian variety with a $\mu_m$-action $\tau_2$.
On $$
\bigl(\p^{m-1}\times A, \Delta:=D_0\times A+\dots + D_{m-1}\times A\bigr)
$$
we have a $\mu_m$-action  generated by $\tau:=(\tau_1, \tau_2)$.

Let $X_1:=\bigl(\p^{m-1}\times A\bigr)/\langle\tau\rangle$.
The quotient of the boundary  $\Delta$ has only 1 component
but it has complicated self-intersections, hence it is not dlt.
Let $(X,\Delta_X)$ be a dlt model.

We see that the minimal lc centers are isomorphic to $(A, 0)$
and the different $\p^1$-linkings between them
differ from each other by a power of $\tau_2$. 
\end{exmp}

\section{Poincar\'e residue map}{\ }

\begin{defn}\label{pres.defn.3.3.say}
Let $(X,\Delta)$ be a dlt pair and
$Z\subset X$ an lc center. As in Definition \ref{doff.defn}, if
$\omega_X^{[m]}(m\Delta)$ is locally free, then, 
by iterating the usual
Poincar\'e residue maps for divisors, we get
a {\it Poincar\'e residue map}
$$
\res_{X\to Z}^m: \omega_X^{[m]}(m\Delta)|_Z\stackrel{\cong}{\longrightarrow} 
\omega_{Z}^{[m]}(m\cdot \diffg_Z\Delta).
\eqno{(\ref{pres.defn.3.3.say}.1)}
$$
(This is well defined if $m$ is even, defined only
up to  sign if $m$ is odd.)

Let $f:(X,\Delta)\to Y$ be a dlt, weak   crepant log structure.
Choose $m>0$ even such that
$\omega_X^{[m]}(m\Delta)\sim f^*L$ for some line bundle $L$ on $Y$. 
Let $Z\subset X$ be an lc center of $(X,\Delta)$. 
We can view the  Poincar\'e residue map  as
$$
\res_{X\to Z}^m: f^*L|_Z\cong 
\omega_X^{[m]}(m\Delta)|_Z\stackrel{\cong}{\longrightarrow} 
\omega_{Z}^{[m]}(m\cdot \diffg_Z\Delta).
\eqno{(\ref{pres.defn.3.3.say}.2)}
$$
\end{defn}

The following result shows, that, for minimal lc centers,
(\ref{pres.defn.3.3.say}.2) is essentially
independent of the choice of $Z$.

\begin{prop}\label{p.res.mincent.prop}
Let $f:(X,\Delta)\to Y$ be a dlt  crepant log structure.
 Choose $m>0$ even such that
$\omega_X^{[m]}(m\Delta)\cong f^*L$ for some line bundle $L$ on $Y$.  
Let $Z_1, Z_2$ be minimal  lc centers  of  $(X,\Delta)$ 
such that  $f(Z_1)=f(Z_2)$.
 Then there is a birational map
$\phi:Z_{2}\map Z_{1}$ such that the following diagram commutes
$$
\begin{array}{rcl}
\omega_X^{[m]}(m\Delta) & \cong f^*L \cong &\omega_X^{[m]}(m\Delta) \\[1ex]
{\res_{X\to Z_{1}}^m}\downarrow {\quad} && {\quad} 
\downarrow{\res_{X\to Z_{2}}^m}\\[1ex]
\omega_{Z_{1}}^{[m]}\bigl(m\diffg_{Z_1}\Delta\bigr) & 
\stackrel{\phi^*}{\longrightarrow}&
 \omega_{Z_{2}}^{[m]}\bigl(m\diffg_{Z_2}\Delta\bigr)
\end{array}
\eqno{(\ref{p.res.mincent.prop}.1)}
$$
\end{prop}

Proof. By Theorem \ref{chains.of.lc.centers.thm}
it is sufficient to prove this in case
there is an 
lc center $W$ that  is birational to a $\p^1$-bundle $\p^1\times U$
with  $Z_1, Z_2$ as sections.
Thus projection to $U$ provides a birational isomorphism
$\phi:Z_2\map Z_1$.

Since  $\res^m_{X\to Z_i}=\res^m_{W\to Z_i}\circ\res^m_{X\to W}$, 
we may assume that $X=W$.
The sheaves in (\ref{p.res.mincent.prop}.1) are torsion free,
hence it is enough to check commutativity after localizing at
the generic point of $U$.
This reduces us to the case when $W=\p^1_L$ with coordinates
$(x{:}y)$, $Z_1=(0{:}1)$ and
$Z_2=(1{:}0)$. A generator of
$H^0\bigl(\p^1, \omega_{\p^1}(Z_1+Z_2)\bigr)$ is $dx/x$
which has residue $1$ at $Z_1$ and $-1$ at $Z_2$. 
Thus (\ref{p.res.mincent.prop}.1) commutes for $m$ even and anti-commutes
for $m$ odd.
 \qed

\begin{rem}\label{p.res.mincent.prop.rem}
 By Proposition \ref{p.res.mincent.prop} we get a 
Poincar\'e residue map as stated in
(\ref{main.intro.new.thm}.4) but it is not yet completely canonical.
We think of $(Z,\Delta_Z)$ as an element of a 
crepant birational equivalence
class, thus so far $\res^m$ is defined only up to the action of
$\bir^c_Y (Z,\Delta_Z)$. 
However, by Theorem \ref{dlt.source.defn.3}, the image of this action is
a finite group of $r$th roots of unity for some $r$.
Thus the $\bir^c_Y (Z,\Delta_Z)$ action is trivial on 
$\omega_Z^{[mr]}(mr\Delta_Z)$ hence
$$
\res^{mr}:  \omega_X^{[mr]}(mr\Delta)|_Z \cong \omega_Z^{[mr]}(mr\Delta_Z)
\eqno{(\ref{p.res.mincent.prop.rem}.1)}
$$
is completely canonical. 
Assume next that $\omega_X^{[mr]}(mr\Delta)\sim f^*L$.
Let us factor $f|_Z:Z\to f(Z)$ using $g:Z\to W$ and 
 the normalization  $n:W\to f(Z)$.
 Then
we can push forward (\ref{p.res.mincent.prop.rem}.1) to get an
isomorphism
$$
 n^*L\cong
\bigl(g_*\omega_Z^{[m]}(m\Delta_Z)\bigr)^{\rm inv}
\eqno{(\ref{p.res.mincent.prop.rem}.2)}
$$
where  the exponent  {\rm inv} denotes the invariants under the
action of the group of birational self-maps $\bir_Y(Z,\Delta_Z)$.
This shows the second isomorphism in
(\ref{main.intro.new.thm}.4).
\end{rem}

\begin{notation} \label{homothety.adj.defn}
Let $(Y,\Delta_Y)$ be lc
and  $(X,\Delta_X)\to (Y,\Delta_Y)$  a crepant, dlt model.
Let $W\subset Y$ be an lc center of $(Y,\Delta_Y)$ and $Z\subset X$  minimal 
(with respect to inclusion) among the lc centers of
$(X,\Delta_X)$ that dominate $W$.
By  Definition \ref{pres.defn.3.3.say}, we obtain a
Poincar\'e residue map  $\res_{X\to Z}$.

Let $D\subset\rdown{\Delta_Y}$ be a divisor
with normalization $\pi:D^n\to D$.
Let $D_X\subset X$ be its birational transform on $X$
and set $\Delta_{D_X}:=\diffg_{D_X}\Delta_X$. 
Let  $W_D\subset D^n$ be  
an lc center of  $\bigl(D^n, \diffg_{D^n}\Delta_Y\bigr)$.
Then $W_X:=\pi(W_D)$ is  an lc center of  $\bigl(Y,\Delta_Y\bigr)$.
Choose  minimal   lc centers $Z_X\subset X$ (resp.\ $Z_D\subset D_X$) 
dominating $W_X$ (resp.\ $W_D$).
\end{notation}

\begin{thm} \label{homothety.adj.thm}
Notation and assumptions as above. Then
 there is a birational map
$\phi:Z_D\map Z_X$ such that for $m$ sufficiently divisible, 
the following diagram commutes
$$
\begin{array}{ccc}
\omega_{X}^{[m]}(m\Delta_X)  & 
\stackrel{\res^m_{X\to D_X}}{\longrightarrow} &
\omega_{D_X}^{[m]}(m\Delta_{D_X})
\\[1ex]
{\res_{X\to Z_X}^m}\downarrow {\quad} && 
{\quad} \downarrow{\res_{ D_X\to Z_D}^m}\\[1ex]
 \omega_{Z_X}^{[m]}\bigl(m\diffg_{Z_X}\Delta_X\bigr) & 
\stackrel{\phi^*}{\longrightarrow}&
\omega_{Z_D}^{[m]}\bigl(m\diffg_{Z_D}\Delta_{D_X}\bigr)
\end{array}
$$
\end{thm}

Proof.
If we choose $Z_X$ as the image of $Z_D$, this holds by the
definition of the higher codimension residue maps.
This and  Proposition \ref{p.res.mincent.prop}
proves the claim for every other choice of $Z_X$.
\qed

\section{Sources and Springs}

\begin{defn}\label{dlt.source.defn}
Let $f:(X,\Delta)\to S$ be a crepant, dlt log structure
and $Z\subset S$ an lc center.
 An lc center
 $Z'$  of  $(X,\Delta)$ is called a {\it   source of $Z$}
if $f(Z')=Z$ and $Z'$ is minimal
(with respect to inclusion) among the lc centers  that
dominate $Z$.

By restriction we have
$f|_{Z'}:\bigl(Z', \diffg_{Z'}\Delta\bigr)\to Z$
and $K_{Z'}+\diffg_{Z'}\Delta\sim_{f,\q} 0$.
By adjunction, there is a one-to-one correspondence
between  lc centers of 
 $\bigl(Z', \diffg_{Z'}\Delta\bigr)$ 
and  lc centers of  $(X,\Delta)$ that are contained in $Z'$.
Thus $Z'$ is a  source of $Z$ iff
the general fiber of 
$\bigl(Z', \diffg_{Z'}\Delta\bigr)\to Z$
 is klt.

By Theorem \ref{chains.of.lc.centers.thm} all  sources of $Z$ are birational to
 each other
(as weak crepant log structures over $Z$). 
Any representative of their birational equivalence class
will be denoted by $\src(Z,X,\Delta)$.
One can choose a representative $(S^t,\Delta^t)\to Z$
whose generic fiber is terminal. Such  models are still not unique,
but their generic fibers are isomorphic outside codimension 2 sets.
However, if there is an irreducible component of $\Delta^t$
whose coefficient is 1 (these can not dominate $Z$) then it does not seem
possible to choose a sensible subclass of models that are 
 isomorphic to each other outside codimension 2 sets.

Note further that by Remark \ref{mincent.bir.inv}, if two crepant log structures
 $f_i:(X_i,\Delta_i)\to Y$ are crepant birational over $Y$, then
 $\src(Z,X_1,\Delta_1)$ is crepant birational to  $\src(Z,X_2,\Delta_2)$.

One can uniquely factor $f|_{Z'}$ as
$$
f|_{Z'}: \bigl(Z', \diffg_{Z'}\Delta'\bigr)=\src(Z,X,\Delta)
\stackrel{c_Z}{\longrightarrow}
\tilde Z' \stackrel{p_Z}{\longrightarrow} Z
\eqno{(\ref{dlt.source.defn}.1)}
$$
where $\tilde Z'$ is normal, $p_Z$ is finite and
$c_Z$ has connected fibers.

Thus in (\ref{dlt.source.defn}.1), $\tilde Z'$ is uniquely defined
up to isomorphism over $Z$. Any representative of its isomorphism class
will be denoted by $\spr(Z,X,\Delta)$ and called the
{\it spring} of $Z$.

Define the group of {\it source-automorphisms} of $\spr(Z,X,\Delta)$ as 
$$
\aut^s \spr(Z,X,\Delta):=
\im\bigl[\bir^c_k\src( Z,X, \Delta)\to \aut_k\spr(Z,X,\Delta) \bigr].
\eqno{(\ref{dlt.source.defn}.2)}
$$
By Theorem \ref{dlt.source.defn.3}, if
 $K_X+\Delta$ is ample then $\aut^s \spr(Z,X,\Delta)$ is finite
 for every lc center $Z\subset X$.

Let $(Y,\Delta)$ be lc and  $f:(X,\Delta_X)\to (Y,\Delta)$
a dlt model. Let $Z\subset Y$ be an lc center of 
 $(Y,\Delta)$.
As noted above,  the source $\src(Z,X,\Delta_X)$ of $Z$  depends 
only on  $(Y,\Delta)$ but not on the choice of $(X,\Delta_X)$.
Thus we also use $\src(Z,Y,\Delta)$ (resp.\ $\spr(Z,Y,\Delta)$)
to denote the source (resp.\ spring) of $Z$. 
\end{defn}

Next we prove (\ref{main.intro.new.thm}.5).

\begin{prop}\label{min.dlt.source.thm}
Let $f:(X,\Delta)\to Y$ be a crepant log structure
 and $Z\subset Y$ an lc center.
  Then the  field extension  $k\bigl(\spr(Z,X,\Delta)\bigr)/k(Z)$ is Galois and
$$
\gal\bigl(\spr(Z,X,\Delta)/Z\bigr)\subset \aut^s \spr(Z,X,\Delta).
$$
\end{prop}

Proof. We may localize at the generic
point of $Z$. Thus we may assume that $Z$ is a point 
 and then prove the following more precise result.

\begin{lem} \label{min.ceneters.sdlt.galois.lem}
Let $g:(X,\Delta)\to Y$ be a weak crepant log structure  over a field $k$.
Assume that $(X,\Delta)$ is dlt and $X$ is $\q$-factorial.
Let  $z\in Y$  be an lc center such that $g^{-1}(z)$ is connected
(as a $k(z)$-scheme).
 Then there is a unique smallest
 finite field extension $K(z)\supset k(z)$ such that
the following hold.
\begin{enumerate}
\item Every  lc center of  $(X_{\bar k},\Delta_{\bar k})$
that intersects $g^{-1}(z)$ 
is defined over $K(z)$.
\item 
Let $W_{\bar z}\subset Y_{\bar k}$ be a minimal lc center contained in  
 $g^{-1}(z)$.  Then
$K(z)=k_{ch}(W_{\bar z})$, the field of definition of $W_{\bar z}$. 
\item $K(z)\supset k(z)$ is a Galois extension.
\item 
Let $W_{z}$ be a minimal lc center contained in $g^{-1}(z)$.
Then 
$$
\bir^c_{k(z)}\bigl(W_z, \diffg_{W_z}\Delta\bigr)
\to\gal\bigl(K(z))/ k(z)\bigr)
\qtq{is surjective.}
$$
\end{enumerate}
\end{lem}

Proof. There are only finitely many lc centers 
and a conjugate of an lc center is also an lc center.
Thus the  field of definition of any lc center 
 is  a finite  extension of $k$. Since
$K(z)$ is the composite of some of them,
 it is  finite over $k(z)$.

 Let $W_{\bar z}\subset X_{\bar k}$ be a minimal lc center  contained in $g^{-1}(z)$
and $k_{ch}(W_{\bar z})$ its
field of definition. Let $D_i\subset \rdown{\Delta}$ be the
irreducible components that contain $W_{\bar z}$. 
Each $D_i$ is smooth at the generic point of $W_{\bar z}$, hence
the $\bar k$-irreducible component of $D_i$ that contains $W_{\bar z}$
is also defined over $k_{ch}(W_{\bar z})$. Thus every  lc center of
$(X_{\bar k},\Delta_{\bar k})$ containing $W_{\bar z}$ is also
defined over $k_{ch}(W_{\bar z})$. Therefore, any lc center that is
$\p^1$-linked to $W_{\bar z}$ is defined over $k_{ch}(W_{\bar z})$. 
By Theorem \ref{chains.of.lc.centers.thm} this implies that
 every  lc center of  $(X_{\bar k},\Delta_{\bar k})$
that intersects $g^{-1}(z)$ is defined over $k_{ch}(W_{\bar z})$, hence
 $k_{ch}(W_{\bar z})\supset K(z)$.
By construction, $k_{ch}(W_{\bar z})\subset  K(z)$, thus
$k_{ch}(W_{\bar z})=  K(z)$.

A conjugate of $W_{\bar z}$ over $k(z)$ is defined over the corresponding 
conjugate field of $k_{ch}(W_{\bar z})$. By the above, every
conjugate of  the field of $k_{ch}(W_{\bar z})$  over $k(z)$ is itself,
hence $k_{ch}(W_{\bar z})=  K(z)$ is Galois over $k(z)$.

Finally, in order to see (4), fix 
$\sigma\in \gal\bigl(K(z)/ k(z)\bigr)$
and let $W_{\bar z}^{\sigma}$ be the corresponding conjugate of $W_{\bar z}$.
By Theorem \ref{chains.of.lc.centers.thm}, $W_{\bar z}^{\sigma}$ and  $W_{\bar z}$
are $\p^1$-linked over $K(z)$; fix one such
$\p^1$-link.
The union of the conjugates of this $\p^1$-link over $k(z)$
 define an element of $\bir^c_{k(z)}\bigl(W_z, \diffg_{W_z}\Delta\bigr)$
which induces $\sigma$ on $K(z)/ k(z)$.
(The $\p^1$-link is not unique, hence the lift is not unique.
Thus in (4) we only claim surjectivity, not a splitting.)
\qed
\medskip

We also note the following direct consequence of 
Corollary \ref{lc.cnet.crep.adj.cor}.

\begin{cor}[Adjunction for sources]\label{adjunction.for.springs}
Let $(X,D+\Delta)$ be lc and  $n:D^n\to D$ the normalization.
Let $Z_D\subset D^n$ be an lc center of
$\bigl(D^n, \diff_{D^n}\Delta\bigr)$ and $Z_X:=n(Z_D)$ its image in $X$. Then
\begin{enumerate}
\item $\src\bigl(Z_D, D^n, \diff_{D^n}\Delta\bigr)\simcb
\src\bigl(Z_X, X, D+\Delta\bigr)$ and
\item $\spr\bigl(Z_D, D^n, \diff_{D^n}\Delta\bigr)\cong
\spr\bigl(Z_X, X, D+\Delta\bigr)$. \qed
\end{enumerate}
\end{cor}

\section{Applications to slc pairs}

\begin{say}[Normalization of slc pairs]\label{slc.norm.basic.say}
Let $(X,\Delta)$ be a semi log canonical pair.
Let $\pi:\bar X\to X$ denote the normalization of $X$,
$\bar\Delta$ the divisorial part of $\pi^{-1}(\Delta)$
and $\bar D\subset \bar X$ the conductor of $\pi$.
Since $X$ is seminormal,  $\bar D$ is reduced.
 $X$ has an ordinary node at a codimension 1  singular point,
thus interchanging the 
two preimages of the node gives an involution
$\tau$ of the normalization $n:\bar D^n\to \bar D$. 
This gives an injection
$$
\left\{
\begin{array}{c}
\mbox{slc  pairs $(X,\Delta)$}
\end{array}
\right\}
\quad  \into \quad 
\left\{\begin{array}{c}
\mbox{lc  pairs  $\bigl(\bar X, \bar D+\bar\Delta\bigr)$}\\
\mbox{plus an involution $\tau$ of $\bar D^n$}
\end{array}
\right\}.
\eqno{(\ref{slc.norm.basic.say}.1)}
$$
For many purposes, it is important to understand the image of
this map. That is, we would like to know which
quadruples  $\bigl(\bar X, \bar D+\bar\Delta, \tau\bigr)$
correspond to an slc  pair $(X,\Delta)$.
An easy condition to derive is that
$\tau$ is an involution not just of the variety $\bar D^n$
but of the lc pair $\bigl(\bar D^n, \diff_{\bar D^n}\bar\Delta\bigr)$.
Thus we obtain a refined version of the map
$$
\left\{
\begin{array}{c}
\mbox{slc  pairs $(X,\Delta)$}
\end{array}
\right\}
\quad  \into \quad 
\left\{
\begin{array}{c}
\mbox{lc  pairs  $\bigl(\bar X, \bar D+\bar\Delta\bigr)$ plus an }\\
\mbox{involution $\tau$ of
$\bigl(\bar D^n, \diff_{\bar D^n}\bar\Delta\bigr)$}
\end{array}
\right\}.
\eqno{(\ref{slc.norm.basic.say}.2)}
$$
For surfaces, the above constructions are discussed in
\cite[Sec.12]{k-etal}. The higher dimensional generalizations
are straightforward; see \cite[Chap.5]{kk-singbook}.

 There are three  major issues involved in trying to prove
that the map (\ref{slc.norm.basic.say}.2) is surjective.
\medskip

\ref{slc.norm.basic.say}.3.1.  
{\it Does $\tau$ generate a finite  equivalence relation?}  

The normalization $n:\bar D^n\to \bar D\to \bar X$ and
$\tau$ generate an equivalence relation $R(\tau)$,
called the {\it gluing relation},  on the points of
$\bar X$ by declaring $n(p)\sim n(\tau(p))$ for every $p\in \bar D^n$.
It is easy to see (cf.\ \cite{k-q}) that
$R(\tau)$ is a set-theoretic, pro-finite, algebraic equivalence relation.
That is, one can give $R(\tau)$ by
 countably many subschemes  
$$
\bigl\{R_i\subset \bar X\times \bar X : i\in I\}
$$
such that 
$\cup_i R_i(K)\subset \bar X(K)\times \bar X(K)$
is an equivalence relation on  $\bar X(K)$ for every algebraically closed field
$K$
and the coordinate projections induce finite morphisms
$$
\pi_1: R_i\to \bar X\qtq{and} \pi_2 :R_i\to \bar X.
$$
(One can make the $R_i$ unique if we choose them irreducible, reduced
and assume that none of them  contains another.)

It is clear that if $X$ exists then every 
 equivalence class  of $R(\tau)$ is contained in a
fiber of $\pi:\bar X\to X$. In particular, if $X$ exists
then the $R(\tau)$-equivalence classes are finite.
Equivalently, $I$ is a finite set.

In general 
the $R(\tau)$-equivalence classes need not be finite.
Moreover, non-finiteness can appear in  high codimension.
This is the question that we will study here using
the sources of lc centers, especially their
Galois property  (\ref{main.intro.new.thm}.5).

A closely related example is given by \cite{Bogomolov-Tschinkel}:
 there is a smooth curve $D$ of genus
$\geq 2$ and a finite relation  $ R_0\subset D\times D$
such that both projections $R_0\rightrightarrows D$ are \'etale
yet $R_0$ generates a non-finite  equivalence relation.
\medskip

\ref{slc.norm.basic.say}.3.2.  
{\it Constructing  $(X,\Delta)$ from 
$\bigl(\bar X, \bar D+\bar\Delta, \tau\bigr)$.}  

Following the method of \cite{k-q},
  it is proved in \cite[Chap.5]{kk-singbook},
that if the $R(\tau)$-equivalence classes are finite,
then $(X,\Delta)$ exists.
\medskip

\ref{slc.norm.basic.say}.3.3. {\it Is $K_X+\Delta$ a  $\q$-Cartier divisor?} 

The answer turns out to be yes, see  \cite[Chap.5]{kk-singbook},
but my proof, using Poincar\'e residue maps and 
Theorem \ref{dlt.source.defn.3}, 
 is somewhat indirect. 
\end{say}

As a consequence we
 obtain that (\ref{slc.norm.basic.say}.2) is
 one-to-one for pairs with ample log canonical class.

\begin{thm} \label{char.of.slcmod.from.nomr}
Taking the normalization gives a one-to-one correspondence between
the following two sets, where $X,\bar X$ are projective schemes over a
field. 
$$
\left\{
\begin{array}{c}
\mbox{slc  pairs $(X,\Delta)$}\\
\mbox{such that}\\
\mbox{$K_X+\Delta$  is ample}
\end{array}
\right\}
\quad \cong\quad
\left\{
\begin{array}{c}
\mbox{lc  pairs  $\bigl(\bar X, \bar D+\bar\Delta\bigr)$ such that  }\\
\mbox{ $K_{\bar X}+ \bar D+\bar\Delta$  is ample plus an }\\
\mbox{involution $\tau$ of
$\bigl(\bar D^n, \diff_{\bar D^n}\bar\Delta\bigr)$}
\end{array}
\right\}.
$$
\end{thm}

This can be extended to the relative case as follows.

\begin{thm} \label{char.of.slcmod.from.nomr.rel}
Let $S$ be a scheme which is essentially of finite type over a
field.
Taking the normalization gives a one-to-one correspondence between
the following two sets.
\begin{enumerate}
\item Slc  pairs $(X,\Delta)$ such that $X/S$ is proper and
$K_X+\Delta$  is ample on the generic fiber of $W\to S$
for every lc center $W\subset X$.
\item Lc pairs $\bigl(\bar X, \bar D+\bar\Delta\bigr)$ 
such that $\bar X/S$ is proper and
$K_{\bar X}+ \bar D+\bar\Delta$  is ample on the generic fiber of $\bar W\to S$
for every lc center $\bar W\subset \bar X$, plus an involution $\tau$ of
$\bigl(\bar D^n, \diff_{\bar D^n}\bar\Delta\bigr)$.
\end{enumerate}
Furthermore, the cases when $K_X+\Delta$  is ample on $X/S$
correspond to the cases when $K_{\bar X}+ \bar D+\bar\Delta$  is ample on
$\bar X/S$.
\end{thm}

As we noted in (\ref{slc.norm.basic.say}.3), the following
result implies  Theorem \ref{char.of.slcmod.from.nomr}. 

\begin{prop}\label{char.of.slcmod.from.prop}
Let $\bigl(\bar X, \bar D+\bar\Delta\bigr)$ be an lc pair
and $\tau$  an involution  of
$\bigl(\bar D^n, \diff_{\bar D^n}\bar\Delta\bigr)$. 

Assume that $X$ is proper over a base scheme $S$ that
 is essentially of finite type over a
field. Assume furthermore that
$K_{\bar X}+ \bar D+\bar\Delta$  is ample on the generic fiber of $\bar W\to S$
for every lc center $\bar W\subset \bar X$.

Then the gluing relation $R(\tau)$,
defined  in (\ref{slc.norm.basic.say}.3.1),  is finite.
\end{prop}

This in turn will be derived from Theorem
\ref{group.on.spring.thm} on the gluing relation $R(\tau)$ which applies
whether $K_{\bar X}+\bar D+\bar\Delta$ is ample or not.

\begin{defn} Let $Y$ be a normal scheme and
$R=\cup_{i\in I} R_i\subset Y\times Y$ a
set-theoretic, pro-finite, algebraic equivalence relation
where the $R_i$ are irreducible.

 $R$ is called a {\it groupoid} if every $R_i$ is the graph of
an isomorphism between two irreducible components of $Y$.

Let $Y^j\subset Y$ be an irreducible component.
The   restriction of $R$ to $Y^j$ is
$R^j:=R\cap  \bigl(Y^j\times Y^j\bigr)$.  If $R$ is a groupoid
then one can identify $R^j$ with a subgroup
of $\aut(Y^j)$ called the {\it stabilizer} of $Y^j$ in $R$.
\end{defn}

We are now ready to formulate and prove a
structure theorem for  gluing relations.
Roughly speaking, we prove that for every lc center
$\bar W\subset \bar X$ there is a ``canonically'' defined
finite cover  $p:\tilde W\to \bar W$ such that
$(p\times p)^{-1}\bigl(R(\tau)\cap (\bar W\times \bar W)\bigr)$
is a  groupoid and the stabilizer 
action  $W$  is compatible with $p^*\bigl(K_{\bar X}+\bar D+\bar\Delta\bigr)$.
The compatibility condition 
is somewhat delicate to state. Thus I give the actual construction
of $\tilde W$ and then specify the compatibility condition
for that particular case.


\begin{notation}\label{group.on.spring.not}
 Let  $(X, \Delta)$ be lc.
Let $S_i^* (X, \Delta)$ be the union of all $\leq i$-dimensional lc centers of
$\bigl( X,  \Delta\bigr)$ and set
$S_i (X, \Delta):=S_i^*(X, \Delta)\setminus S_{i-1}^*(X, \Delta)$.
Let $Z_{ij}^0\subset S_i(X, \Delta)$ be the irreducible components.
The closure $Z_{ij}$ of  $Z_{ij}^0$
 is an lc center of $\bigl( X,  \Delta\bigr)$,
hence it has a spring  $p_{ij}:\spr(Z_{ij},X,  \Delta)\to Z_{ij}$.
Set 
$\spr(Z^0_{ij},X,  \Delta):=p_{ij}^{-1}Z_{ij}^0$ and
$$
\spr_i\bigl(X,  \Delta\bigr):=\amalg_j  \spr(Z^0_{ij},X,  \Delta).
$$
Let $p_i:\spr_i\bigl(X,  \Delta\bigr)\to S_i(X, \Delta)$ 
be the induced morphism. Then $p_i$ is finite, surjective and
universally open since $S_i(X, \Delta)$ is normal.
Furthermore, $p_i$ is Galois over every  $Z_{ij}$  by 
Proposition \ref{min.dlt.source.thm}.
\end{notation}

\begin{thm}\label{group.on.spring.thm}
 Let  $\bigl(X, D+\Delta\bigr)$ be lc,
 $\tau$  an involution  of
$\bigl(D^n, \diff_{ D^n}\Delta\bigr)$ and 
 $ R(\tau)\subset  X\times  X$
 the corresponding   equivalence relation as in 
(\ref{slc.norm.basic.say}.3.1).
Let $p_i:\spr_i\bigl(X,  D+\Delta\bigr)\to S_i$  be as above. Then 
\begin{enumerate}
\item $(p_i\times p_i)^{-1}\bigl(R(\tau)\cap 
(S_i(X, \Delta)\times S_i(X, \Delta))\bigr)$
is a groupoid on $\spr_i\bigl(X,  D+\Delta\bigr)$.
\item For every  irreducible component
$Z_{ij}^0\subset S_i(X, \Delta)$,
the stabilizer of  its spring $\spr(Z^0_{ij},X,  D+\Delta)
\subset \spr_i\bigl(X,  D+\Delta\bigr)$
is a subgroup of the source-automor\-phism group 
$\aut^s \spr(Z_{ij},X,D+\Delta)$.
\end{enumerate}
\end{thm}

Proof.  We need to describe how the generators of $R(\tau)$
pull back to the spring $\spr_i\bigl(X,  D+\Delta\bigr)$.

First, the preimage of the diagonal of $Z_{ij}^0\times Z_{ij}^0$
is a group $\Gamma(G_{ij})$ and 
$G_{ij}=\gal\bigl(\spr(Z_{ij},X,  D+\Delta)/Z_{ij}\bigr)$
  is a subgroup of
 $\aut^s \spr(Z_{ij},X,D+\Delta)$
by Proposition \ref{min.dlt.source.thm}.

Second, let  $Z_{ijk}\subset D^n$ be an irreducible component of 
the preimage of $Z_{ij}$. Then  $Z_{ijk}$ is an lc center
of $\bigl(D^n,  \diff_{ D^n}\Delta\bigr)$ and  
$$
\src\bigl(Z_{ijk}, D^n,  \diff_{ D^n}\Delta\bigr)\simcb 
\src\bigl(Z_{ij}, X,  D+\Delta\bigr)
$$
by  Corollary \ref{adjunction.for.springs}.
Thus, for each ${ijk}$,
 the isomorphism $\tau:D^n\cong D^n$ lifts to
isomorphisms 
$$
\tau_{ijkl}:\spr(Z^0_{ij},X,  D+\Delta)\cong \spr(Z^0_{il},X,  D+\Delta).
$$
 Given $ijk$, the value
of $l$ is determined 
 by  $Z_{il}:=n\bigl(\tau(Z_{ijk})\bigr)$, 
but the lifting is defined
only up to left and right multiplication by elements of 
$G_{ij}$ and $G_{il}$. 

Thus 
$(p_i\times p_i)^{-1}\bigl(R(\tau)\cap
 (S_i(X, \Delta)\times S_i(X, \Delta))\bigr)$
is the groupoid generated by the $G_{ij}$
and the $\tau_{ijkl}$, hence
the stabilizer of $\spr(Z^0_{ij},X,  D+\Delta)$
is  generated by the groups
$\tau_{ijkl}^{-1}G_{il}\tau_{ijkl}$. The latter  are all
subgroups of $\aut^s \spr(Z_{ij},X,D+\Delta)$.\qed

\begin{say}[Proof of Proposition \ref{char.of.slcmod.from.prop}]
Since $\spr_i(X,  D+\Delta)$ has finitely many
irreducible components, the groupoid is finite iff the
stabilizer of each $\spr(Z^0_{ij},X,  D+\Delta)$ is finite.
By  Theorem \ref{group.on.spring.thm} this holds if
the groups $\aut^s \spr(Z_{ij},X,D+\Delta)$   are finite.

The automorphism group of a variety $\tilde Z$ over
a base scheme $S$ injects into the automorphism group of
the generic fiber $\tilde Z_{gen}$.

By assumption,  $K_{ \bar X}+ \bar D+\bar\Delta$ is ample on 
the generic fiber of $Z_{ij}\to S$,  thus  Theorem \ref{dlt.source.defn.3}
 implies that  each $\aut^s \spr(Z_{ij},X,D+\Delta)$ 
is finite. \qed
\end{say}

 \begin{ack} This paper was written while I visited  RIMS, Kyoto University.
I thank  S.~Mori and S.~Mukai for the invitation and their hospitality.
I am grateful to O.~Fujino, 
C.~Hacon, S.~Kov\'acs, S.~Mori, Y.~Odaka, V.~Shokurov and C.~Xu for
many   comments and corrections.
Partial financial support  was provided by  the NSF under grant number 
DMS-0758275.
\end{ack}


\begin{thebibliography}{BCHM10}

\bibitem[BCHM10]{bchm}
Caucher Birkar, Paolo Cascini, Christopher~D. Hacon, and James McKernan,
  \emph{Existence of minimal models for varieties of log general type}, J.
  Amer. Math. Soc. \textbf{23} (2010), no.~2, 405--468. \MR{2601039}

\bibitem[BPV84]{bpv}
W.~Barth, C.~Peters, and A.~V{an de Ven}, \emph{Compact complex surfaces},
  Ergebnisse der Mathematik und ihrer Grenzgebiete (3), vol.~4,
  Springer-Verlag, Berlin, 1984. \MR{749574 (86c:32026)}

\bibitem[BT09]{Bogomolov-Tschinkel}
Fedor Bogomolov and Yuri Tschinkel, \emph{Co-fibered products of algebraic
  curves}, arXiv.org:0902.0534, 2009.

\bibitem[FG10]{fuj-gon-cbf}
Osamu {Fujino} and Y.~{Gongyo}, \emph{{Log pluricanonical representations and
  abundance conjecture}}, ArXiv e-prints (2010).

\bibitem[Fuj07]{fuj-in-corti}
Osamu Fujino, \emph{What is log terminal?}, Flips for 3-folds and 4-folds,
  Oxford Lecture Ser. Math. Appl., vol.~35, Oxford Univ. Press, Oxford, 2007,
  pp.~49--62. \MR{2359341}

\bibitem[Fuj10]{fujino-ssmmp}
\bysame, \emph{{Semi-stable minimal model program for varieties with trivial
  canonical divisor}}, ArXiv e-prints (2010).

\bibitem[{Gon}10]{gon-ab}
Y.~{Gongyo}, \emph{{Abundance theorem for numerically trivial log canonical
  divisors of semi-log canonical pairs}}, ArXiv e-prints (2010).

\bibitem[Kaw97]{kaw1}
Yujiro Kawamata, \emph{Subadjunction of log canonical divisors for a subvariety
  of codimension {$2$}}, Birational algebraic geometry ({B}altimore, {MD},
  1996), Contemp. Math., vol. 207, Amer. Math. Soc., Providence, RI, 1997,
  pp.~79--88. \MR{1462926 (99a:14024)}

\bibitem[Kaw98]{kaw-adj}
\bysame, \emph{Subadjunction of log canonical divisors. {II}}, Amer. J. Math.
  \textbf{120} (1998), no.~5, 893--899. \MR{1646046 (2000d:14020)}

\bibitem[KK10]{k-db}
J{\'a}nos Koll{\'a}r and S{\'a}ndor~J. Kov{\'a}cs, \emph{Log canonical
  singularities are {D}u {B}ois}, J. Amer. Math. Soc. \textbf{23} (2010),
  no.~3, 791--813. \MR{2629988}

\bibitem[KM98]{km-book}
J{\'a}nos Koll{\'a}r and Shigefumi Mori, \emph{Birational geometry of algebraic
  varieties}, Cambridge Tracts in Mathematics, vol. 134, Cambridge University
  Press, Cambridge, 1998, With the collaboration of C. H. Clemens and A. Corti,
  Translated from the 1998 Japanese original. \MR{1658959 (2000b:14018)}

\bibitem[Kol86]{k-hdi1}
J{\'a}nos Koll{\'a}r, \emph{Higher direct images of dualizing sheaves. {I}},
  Ann. of Math. (2) \textbf{123} (1986), no.~1, 11--42. \MR{825838 (87c:14038)}

\bibitem[Kol92]{k-etal}
J{\'a}nos Koll{\'a}r (ed.), \emph{Flips and abundance for algebraic
  threefolds}, Soci\'et\'e Math\'ematique de France, 1992, Papers from the
  Second Summer Seminar on Algebraic Geometry held at the University of Utah,
  Salt Lake City, Utah, August 1991, Ast\'erisque No. 211 (1992).
  \MR{94f:14013}

\bibitem[Kol97]{k-quot}
\bysame, \emph{Quotient spaces modulo algebraic groups}, Ann. of Math. (2)
  \textbf{145} (1997), no.~1, 33--79. \MR{1432036 (97m:14013)}

\bibitem[Kol07]{k-adj}
\bysame, \emph{Kodaira's canonical bundle formula and adjunction}, Flips for
  3-folds and 4-folds, Oxford Lecture Ser. Math. Appl., vol.~35, Oxford Univ.
  Press, Oxford, 2007, pp.~134--162. \MR{2359346}

\bibitem[Kol12]{k-q}
\bysame, \emph{Quotients by finite equivalence relations}, Current developments
  in algebraic geometry, Math. Sci. Res. Inst. Publ., vol.~59, Cambridge Univ.
  Press, Cambridge, 2012, With an appendix by Claudiu Raicu, pp.~227--256.
  \MR{2931872}

\bibitem[Kol13]{kk-singbook}
\bysame, \emph{Singularities of the minimal model program}, Cambridge
  University Press, Cambridge, 2013, With the collaboration of S. Kov\'acs (to
  appear).

\bibitem[Mil80]{milne}
James~S. Milne, \emph{\'{E}tale cohomology}, Princeton Mathematical Series,
  vol.~33, Princeton University Press, Princeton, N.J., 1980. \MR{559531
  (81j:14002)}

\bibitem[NU73]{nak-ue}
Iku Nakamura and Kenji Ueno, \emph{An addition formula for {K}odaira dimensions
  of analytic fibre bundles whose fibre are {M}oi\v sezon manifolds}, J. Math.
  Soc. Japan \textbf{25} (1973), 363--371. \MR{0322213 (48 \#575)}

\bibitem[Oda11]{odaka}
Yuji Odaka, \emph{The {GIT} stability of polarized varieties via discrepancy},
  2011.

\bibitem[OX12]{oda-xu}
Yuji Odaka and Chenyang Xu, \emph{Log-canonical modification of singular pairs
  and its applications}, Math. Res. Lett. \textbf{19} (2012), no.~2, 325--334.

\bibitem[Sho13]{shok-source}
V.~V. Shokurov, \emph{Log adjunction: effectiveness and positivity (in
  preparation)}, 2013.

\bibitem[Sza94]{szabo}
Endre Szab{\'o}, \emph{Divisorial log terminal singularities}, J. Math. Sci.
  Univ. Tokyo \textbf{1} (1994), no.~3, 631--639. \MR{1322695 (96f:14019)}

\bibitem[Uen75]{ueno}
Kenji Ueno, \emph{Classification theory of algebraic varieties and compact
  complex spaces}, Lecture Notes in Mathematics, Vol. 439, Springer-Verlag,
  Berlin, 1975, Notes written in collaboration with P. Cherenack. \MR{0506253
  (58 \#22062)}

\end{thebibliography}

\def\cprime{$'$} \def\cprime{$'$} \def\cprime{$'$} \def\cprime{$'$}
  \def\cprime{$'$} \def\cprime{$'$} \def\dbar{\leavevmode\hbox to
  0pt{\hskip.2ex \accent"16\hss}d} \def\cprime{$'$} \def\cprime{$'$}
  \def\polhk#1{\setbox0=\hbox{#1}{\ooalign{\hidewidth
  \lower1.5ex\hbox{`}\hidewidth\crcr\unhbox0}}} \def\cprime{$'$}
  \def\cprime{$'$} \def\cprime{$'$} \def\cprime{$'$}
  \def\polhk#1{\setbox0=\hbox{#1}{\ooalign{\hidewidth
  \lower1.5ex\hbox{`}\hidewidth\crcr\unhbox0}}} \def\cdprime{$''$}
  \def\cprime{$'$} \def\cprime{$'$} \def\cprime{$'$} \def\cprime{$'$}
\providecommand{\bysame}{\leavevmode\hbox to3em{\hrulefill}\thinspace}
\providecommand{\MR}{\relax\ifhmode\unskip\space\fi MR }
\providecommand{\MRhref}[2]{%
  \href{http://www.ams.org/mathscinet-getitem?mr=#1}{#2}
}
\providecommand{\href}[2]{#2}

\noindent Princeton University, Princeton NJ 08544-1000

\begin{verbatim}kollar@math.princeton.edu\end{verbatim}

\end{document}